\documentclass{llncs}

\usepackage{graphicx}

\usepackage{tikz}

\usepackage{color}

\usepackage{subfigure} % Written by Steven Douglas Cochran
\usepackage{amsmath}   % From the American Mathematical Society
\usepackage{amssymb}

\usepackage{float}

\usepackage{ifpdf}

\newcommand{\Cauchy}{L} % {\mathcal L}
\newcommand{\Momentum}{L}

\newcommand{\IR}{\mathbb{R}}

\begin{document}

\title{Combining the band-limited parameterization and Semi-Lagrangian Runge-Kutta integration for efficient PDE-constrained LDDMM}

% \subtitle{Do you have a subtitle?\\ If so, write it here}

\titlerunning{BL parameterization and SL-RK for PDE-LDDMM}        % if too long for running head

\author{Monica Hernandez}


\institute{Computer Sciences Department \\ Aragon Institute on Engineering Research \\ University of Zaragoza \\ mhg@unizar.es}

\date{Received: This work was finished in January 2019  / Accepted: date}

\maketitle

\begin{abstract}
The family of PDE-constrained LDDMM methods is emerging as a particularly interesting approach for physically meaningful diffeomorphic transformations. 
The original combination of Gauss--Newton--Krylov optimization
and Runge-Kutta integration, 
shows excellent numerical accuracy and fast convergence rate. 
However, its most significant limitation is the huge computational complexity, 
hindering its extensive use in Computational Anatomy applied studies. 
This limitation has been treated independently by the problem formulation in the space of band-limited vector fields and Semi-Lagrangian integration.
The purpose of this work is to combine both in three variants of band-limited PDE-constrained LDDMM 
for further increasing their computational efficiency. 
The accuracy of the resulting methods is evaluated extensively.
For all the variants, the proposed combined approach shows a significant increment of the computational efficiency.
In addition, the variant based on the deformation state equation is positioned consistently as the best performing method 
across all the evaluation frameworks in terms of accuracy and efficiency. 
This work was finished on January 2019 and is currently under review.
% \keywords{physically meaningful diffeomorphic registration \and PDE-constrained LDDMM \and Gauss--Newton--Krylov \and 
% optimal control optimization \and band-limited vector fields \and Semi-Lagrangian Runge-Kutta integration }
\end{abstract}

\newpage

\section{Introduction}

Computational Anatomy is a powerful interdisciplinary field for the analysis of anatomical shape variability~\cite{Miller_04,Miller_09}.
This discipline is based on Sir D'Arcy Thompson's original ideas for explaining the similarity of the anatomical shape
of homologous species using the transformations existing between the anatomical structures~\cite{Thompson_17}.
In Computational Anatomy, shape similarity is measured from the diffeomorphic transformations estimated between the anatomies.
These transformations yield a generative model for the analysis of shape variability.
Diffeomorphisms are computed from the anatomical images using diffeomorphic registration methods~\cite{Younes_Book}.

There exists a vast literature on diffeomorphic registration methods with differences in the transformation characterization, 
regularizers, image similarity metrics, optimization methods, and additional constraints~\cite{Sotiras_13}. 
Although the differentiability and invertibility of the transformations constitute crucial features for Computational Anatomy applications, 
the diffeomorphic constraint does not necessarily guarantee that a transformation computed with a given method is physically meaningful for the clinical domain of interest. 
PDE-constrained Large Deformation Diffeomorphic Metric Mapping (PDE-LDDMM) has become relevant in the last decade
for the computation of transformations under plausible physical models of interest~\cite{Younes_07,Hart_09,Ashburner_11,Vialard_11,Mang_15,Mang_16,Mang_17,Zhang_18,Hernandez_19b}. 

Our work focuses on the family of PDE-constrained LDDMM methods pioneered by Hart et al.~\cite{Hart_09} and leading to the relevant contributions 
in~\cite{Mang_15,Mang_16,Polzin_16,Hernandez_19c,Mang_19}.
In this family of methods, the registration problem is approached from an optimal control perspective, where the different physical models are imposed 
directly using the physical PDEs that are attached to the LDDMM variational 
problem using hard constraints.
The numerical optimization is approached using gradient-descent~\cite{Hart_09,Mang_15,Polzin_16} or second-order optimization in the form 
of inexact reduced Newton-Krylov methods~\cite{Mang_15,Mang_16,Hernandez_19c,Mang_19}. 
In particular, the combination of Gauss--Newton--Krylov for optimization, with sophisticated multi-level preconditioners, 
spectral methods for differentiation, and Runge-Kutta schemes for PDE integration,  
shows excellent numerical accuracy and an extraordinarily fast convergence rate.
However, the most significant limitation of Gauss--Newton--Krylov PDE-constrained LDDMM is the huge computational complexity, which hinders 
the extensive use in Computational Anatomy applied studies.
This computational complexity is due to: 
\begin{enumerate}
\item The formulation of the problem in the spatial domain. 
\item The large time-sampling needed for the stability of Runge-Kutta integration.
\end{enumerate}

\noindent Both issues have been treated independently in the literature yielding to PDE-constrained LDDMM methods with increased 
efficiency and an assumable cost in accuracy loss.

\subsection{Computational complexity due to problem formulation}

The computational complexity due to the formulation of the problem in the spatial domain has been successfully reduced using the band-limited vector field 
parameterization proposed in~\cite{Zhang_15,Zhang_18}.
LDDMM methods, and in particular PDE-constrained LDDMM, involve the action of low-pass filters in the optimization update equations of the velocities.
Therefore, the computation of the high-frequency components of high-resolution velocity fields can be omitted since these computations result equal or nearly 
equal to zero by the action of the low-pass filters.
The band-limited vector field parameterization allows a reduction of the dimensionality of the problem that circumvents the high-frequency computations.

The works in~\cite{Hernandez_19,Hernandez_19c} formulate three different variants of PDE-constrained LDDMM in the space of band-limited vector fields
and perform the computations in the GPU.
Some configurations of these variants have been really successful, greatly outperforming the state of the art methods in terms of computational complexity
while keeping a competitive accuracy.

\subsection{Computational complexity due to PDE integration}

Runge-Kutta methods are explicit techniques. Hence they are only conditionally stable.
This means that the time sampling should be selected enough to preserve the Courant-Friedrichs-Lewy (CFL) condition.
For PDE-constrained LDDMM, the time sampling values that guarantee stability are usually large.
As a result, the time and memory requirements of the problem are considerably increased.
In particular, the memory requirements of PDE-constrained LDDMM is increased to limits that hinder the execution on limited memory
devices such as the GPU.
In addition, one can experience that the time-sampling needed for the non-stationary parameterization is much higher than for the 
stationary parameterization, increasing the complexity of an already not particularly memory efficient configuration.
On the other side, when stability is satisfied, the accuracy of PDE-constrained LDDMM is high~\cite{Mang_15,Mang_16,Hernandez_19,Hernandez_19c}.

Semi-Lagrangian methods are semi-implicit techniques that are unconditionally stable.
Therefore, the time sampling can be selected according to accuracy rather than stability considerations.
Semi-Lagrangian methods were originally proposed in the 90's in the context of modelling weather predictions~\cite{Staniforth_91}.
In the context of diffeomorphic registration, the original LDDMM method proposed in~\cite{Beg_05} already used Semi-Lagrangian integration 
for the solution of the transport equation.
The combination of Semi-Lagrangian integration with Runge-Kutta has been recently proposed for solving some time-dependent PDEs.
Runge-Kutta has shown to increase the accuracy of first-order schemes in Semi-Lagrangian integration~\cite{Guo_13}.

The computational complexity in~\cite{Mang_15,Mang_16} due to the use of Runge-Kutta schemes for PDE integration has been successfully reduced 
using Semi-Lagrangian Runge-Kutta integration~\cite{Mang_17b,Mang_19} for the stationary parameterization of diffeomorphisms.
For PDE-constrained LDDMM, the selected time sampling is usually much smaller than the time sampling typically selected with explicit 
schemes, yielding to a considerable reduction of the computational complexity of the problem.
On the other hand, the expected accuracy of PDE-constrained LDDMM is lower than with explicit schemes.
% ~\footnote{This result has not been reported in the state of the art.
% However, it is consistently reported by our extensive evaluation.}.

Beyond the computational complexity improvement through numerical schemes, Mang et al. proposed an efficient implementation of PDE-LDDMM that 
exploits massive CPU based parallel computing architectures~\cite{Mang_16_3D}. The source code has been recently released with~\cite{Mang_19}. 

\subsection{Our contribution}

The purpose of this work is to further increase the computational efficiency of BL PDE-constrained LDDMM by combining the two independent methodological 
approaches of circumventing the huge computational complexity of PDE-constrained LDDMM and to extensively analyze the accuracy of the resulting methods.
We have implemented the band-limited methods in~\cite{Hernandez_19,Hernandez_19c} with the Semi-Lagrangian Runge-Kutta integration scheme originally 
proposed in~\cite{Mang_17b} for the stationary and the non-stationary parameterization of diffeomorphisms.
The resulting methods have been evaluated in five different datasets following the evaluation frameworks in~\cite{Hernandez_19c,Klein_09,Rohlfing_12}.
To our knowledge, this is the first time that Semi-Lagrangian Runge-Kutta integration is implemented in the space of band-limited vector fields.
It is also the first time that Semi-Lagrangian Runge-Kutta integration is used in PDE-LDDMM with the non-stationary parameterization.
Moreover, our work first provides the position achieved by benchmark PDE-constrained LDDMM methods in the ranking of Klein et al. evaluation.
The best performing method of our work coincides with the best performing variant in~\cite{Hernandez_19c}, PDE-constrained LDDMM based on the deformation state equation.
The Semi-Lagrangian Runge-Kutta scheme proposed in this work has shown to outperform the Runge-Kutta scheme in~\cite{Hernandez_19c} in terms of computational efficiency
and accuracy.
Indeed, the best performing PDE-LDDMM variant in this work has recently reached the highest sensitivity (97 \% vs a baseline of 88 \%) in the classification 
of stable vs progressive mild cognitive impaired conversors in the Alzheimer's Disease Neuroimaging Initiative (ADNI) database using convolutional neural networks~\cite{Ramon_20}. 

\subsection{Manuscript organization}

In the following, Section~\ref{sec:PDE-LDDMM} reviews the foundations of PDE-constrained LDDMM, 
with particular emphasis on the band-limited vector field parameterization.
Section~\ref{sec:SL} presents the proposed Semi-Lagrangian Runge-Kutta integration method.
Next, Section~\ref{sec:ExpSetup} details the experimental setup.
Section~\ref{sec:Results} shows the results and Section~\ref{sec:Discussion} discusses the most important highlights. 
Finally, Section~\ref{sec:Conclusions} gathers the most remarkable conclusions of our work.

\section{PDE-constrained LDDMM methods}
\label{sec:PDE-LDDMM}

\subsection{Parameterization in the spatial domain}

Let $\Omega \subseteq \mathbb{R}^d$ be the image domain.
Let $Diff(\Omega)$ be the LDDMM Riemannian manifold of diffeomorphisms and $V$ the tangent space at the identity element.
$Diff(\Omega)$ is a Lie group, and $V$ is the corresponding Lie algebra.
The Riemannian metric of $Diff(\Omega)$ is defined from the scalar product in $V$ 
\begin{equation}
 \langle v, w \rangle_V = \langle \Cauchy v, w \rangle_{L^2} = \int_\Omega \langle \Cauchy v(x), w(x) \rangle d\Omega, 
\end{equation}
where $\Cauchy = (Id - \alpha \Delta)^s, \alpha >0, s \in \mathbb{N}$
is the invertible self-adjoint differential operator associated with the differential structure of 
$Diff(\Omega)$~\cite{Beg_05}.
% We denote with $K$ to the inverse of operator $L$.

Let $I_0$ and $I_1$ be the source and the target images.
PDE-constrained LDDMM is formulated from the minimization of the PDE-constrained variational problem

\begin{equation}
\label{eq:LS-LDDMM}
E(v) = \frac{1}{2} \int_0^1 \langle \Cauchy v_t, v_t \rangle_{L^2} dt + \frac{1}{\sigma^2} \Vert m(1) - I_1\Vert_{L^2}^2, 
\end{equation}
subject to the state equation with state variable $m(t)$
\begin{equation}
\label{eq:StateEquation}
\partial_t m(t) + \nabla m(t) \cdot v_t = 0 \textnormal{ in } \Omega \times (0,1],
\end{equation}
\noindent with initial condition $m(0) = I_0$~\cite{Hart_09,Mang_15}. 

The variational problem is posed in the space of time-varying smooth flows of velocity fields, ${v} \in L^2([0,1],V)$.
Given the smooth flow ${v}:[0,1] \rightarrow V$, $v_t:\Omega \rightarrow \mathbb{R}^{d}$, 
the solution at time $t=1$ to the evolution equation
\begin{equation}
\label{eq:TransportEquation}
\partial_t \phi_t^{v} = v_t \circ \phi_t^{v}  
\end{equation}
with initial condition $\phi_0^{v} = id$ is a diffeomorphism, $\phi_1^{v} \in Diff(\Omega)$.
The transformation $(\phi^{v}_{1})^{-1}$, computed from the minimum of $E({v})$ constrained to the state equation, 
is the diffeomorphism that solves the PDE-constrained LDDMM registration problem between $I_0$ and $I_1$.

The state equation constraint in PDE-constrained LDDMM can be imposed in two more different manners, yielding three
different variants of PDE-constrained LDDMM~\cite{Hernandez_19c}. The second variant is formulated from the 
minimization of Equation~\ref{eq:LS-LDDMM}, where 
\begin{equation}
 m(t) = I_0 \circ \phi(t)
\end{equation}
\noindent and $\phi$ is computed from the deformation state equation
\begin{equation}
\label{eq:DefStateEquation}
\partial_t \phi(t) + D \phi(t) \cdot v_t = 0 \textnormal{ in } \Omega \times (0,1].
\end{equation}
\noindent The third variant is formulated from the minimization of Equation~\ref{eq:LS-LDDMM} subject to the
deformation state equation (Equation~\ref{eq:DefStateEquation}).

The best optimization method for PDE-constrained LDDMM is Gauss--Newton--Krylov~\cite{Mang_15,Hernandez_19c}.
The expressions of the gradient and the Hessian-vector product are derived from the augmented Lagrangian of the 
energy functional and the state equation.
The update equation has the form

\begin{equation}
\label{eq:Newton}
 v^{n+1} = v^n + \epsilon \delta v^n,
\end{equation}
where $\epsilon$ is the update length and $\delta v^n$ is computed from preconditioned conjugate gradient (PCG) on the system
\begin{equation}
\label{eq:KKT}
 H_v E_{aug}( v^n ) \delta  v^n = - \nabla_v E_{aug}(v^n),
\end{equation}
with preconditioner $\Cauchy^{-1}$. 

\subsection{Parameterization in the space of band-limited vector fields}

Let $\widetilde{\Omega}$ be the discrete Fourier domain truncated with frequency bounds $K_1,$ $\dots,$ $K_d$.
We denote with $\widetilde{V}$ the space of discretized band-limited vector fields on $\Omega$ with these frequency bounds. 
The elements in $\widetilde{V}$ are represented in the Fourier domain as $\tilde{v}: \widetilde{\Omega} \rightarrow \mathbb{C}^d$,
$\tilde{v}(k_1, \dots, k_d)$.
% , and in the spatial domain as $\iota(\tilde{v}):\Omega \rightarrow \mathbb{R}^d$,
% \begin{equation}
% \iota(\tilde{v})(x_1, \dots, x_d) = \sum_{k_1 = 0}^{K_1} \dots \sum_{k_d = 0}^{K_d} \tilde{v}(k_1, \dots, k_d) e^{2 \pi i k_1 x_1} \dots e^{2 \pi i k_d x_d}. 
% \end{equation}
\noindent The application $\iota:\widetilde{V} \rightarrow V$ denotes the natural inclusion mapping of $\widetilde{V}$ in $V$.
The application $\pi: V \rightarrow \widetilde{V}$ denotes the projection of $V$ onto $\widetilde{V}$~\cite{Zhang_15,Zhang_18}.

The space $\widetilde{V}$ of band-limited vector fields has a finite-dimensional Lie algebra structure using the truncated 
convolution $\star$ in the definition of the Lie bracket~\cite{Zhang_18}.
We denote with $Diff(\widetilde{\Omega})$ to the finite-dimensional Riemannian manifold of diffeomorphisms on $\widetilde{\Omega}$ 
with corresponding Lie algebra $\widetilde{V}$.
The Riemannian metric in $Diff(\widetilde{\Omega})$ is defined from the scalar product
\begin{equation}
 \langle \tilde{v}, \tilde{w} \rangle_{\tilde{V}} = \langle \tilde{\Cauchy} \tilde{v}, \tilde{w} \rangle_{l^2},  
\end{equation}
\noindent where $\tilde{\Cauchy}$ is the projection of operator $\Cauchy$ in the truncated Fourier domain. 
Similarly, we will denote with $\tilde{\ast}$ to the projection in the truncated Fourier domain of the differential operators $\ast$
involved in the differential equations.
% Similarly, we will denote with $\tilde{K}$, $\widetilde{\nabla}$, and $\widetilde{\nabla \cdot}$ to 
% the projection of operators $K$, $\nabla$, and $\nabla \cdot$ in the truncated Fourier domain.
% In addition, we will denote with $\star$ the truncated convolution.
 
The band-limited PDE-constrained variational problem is given by the minimization of 
\begin{equation}
\label{eq:BLEnergy}
 E(\tilde{v}) = \frac{1}{2} \int_0^1 \langle \tilde{L}\tilde{v}_t, \tilde{v}_t \rangle_{l^2} dt + \frac{1}{\sigma^2} \Vert m(1) - I_1 \Vert_{L^2}^2.
\end{equation}
The first variant is formulated from the minimization of Equation~\ref{eq:BLEnergy} subject to 
\begin{equation}
\partial_t m(t) + \nabla m(t) \cdot \iota(\tilde{v}_t) = 0 \textnormal{ in } \Omega \times (0,1], \\
\end{equation}
\noindent with initial condition $m(0) = I_0$. 
For the second variant, the diffeomorphism is computed from $\phi(t) = id - \iota(\tilde{u})(t)$ where
$\tilde{u}(t)$ is computed from the deformation state equation formulated in displacement field form
\begin{equation}
\label{eq:BLDefStateEquation}
\partial_t \tilde{u}(t) + \widetilde{D} \tilde{u}(t) \star \tilde{v}_t = \tilde{u}(t) \textnormal{ in } \widetilde{\Omega} \times (0,1].
\end{equation}
The third variant is formulated analogously to the spatial case from the minimization of Equation~\ref{eq:BLEnergy} subject to the
deformation state equation (Equation~\ref{eq:BLDefStateEquation})~\cite{Hernandez_19,Hernandez_19c}.

The optimization is approached using Gauss--Newton--Krylov methods in $\widetilde{V}$ with preconditioner $\tilde{\Cauchy}^{-1}$.
The update equation has the form
\begin{equation}
\label{eq:BLNewton}
 \tilde{v}^{n+1} = \tilde{v}^n + \epsilon \delta \tilde{v}^n,
\end{equation}
\noindent where $\delta \tilde{v}^n$ is computed from
\begin{equation}
\label{eq:BLKKT}
 \widetilde{(H_{\tilde{v}} E_{aug}(  \tilde{v}^n))} \delta \tilde{v}^n = - \widetilde{(\nabla_{\tilde{v}} E_{aug}(\tilde{v}^n))}.
\end{equation}

\subsection{BL PDE-constrained LDDMM equations}
\label{sec:OPDE-LDDMM}

% \color{red}
% The state equation constraint in PDE-constrained LDDMM can be imposed in three different manners, yielding three
% different variants of PDE-constrained LDDMM and BL PDE-constrained LDDMM~\cite{Hernandez_19,Hernandez_19c}. 
% In the following, we focus on the equations for the band-limited parameterization.
% \color{black}

\subsubsection{Original BL PDE-constrained LDDMM}

Originally proposed BL PDE-constrained LDDMM uses the state equation in the augmented Lagrangian for the derivation 
of the state and adjoint equations and their incremental counterparts~\cite{Hernandez_19}. 
The gradient and the Gauss-Newton approximation of the Hessian-vector product are given by the equations
\begin{eqnarray}
\widetilde{(\nabla_{\tilde{v}} E_{aug}(\tilde{v}))_t} = \tilde{\Momentum} \tilde{v}_t + \tilde{\lambda}(t) \star \widetilde{\nabla} \tilde{m}(t) \\
\widetilde{ (H_{\tilde{v}} E_{aug}(\tilde{v}))_t } \delta \tilde{v} (t) = \tilde{\Momentum} \delta \tilde{v}(t) + \delta \tilde{\lambda}(t) \star \widetilde{\nabla} \tilde{m}(t),
\end{eqnarray}
\noindent where the projected state variable $\tilde{m}$ and the projected adjoint variable $\tilde{\lambda}$ are computed from $\pi(m)$ and $\pi(\lambda)$, and $m$ and $\lambda$ are computed from
\begin{eqnarray}
\label{eq:BLState}
\partial_t m(t) + \nabla m(t) \cdot \iota(\tilde{v}_t) = 0 \\
\label{eq:BLAdjoint}
-\partial_t \lambda(t) - \nabla \cdot ( \lambda(t) \cdot \iota(\tilde{v}_t) ) = 0.
\end{eqnarray}
\noindent The incremental counterparts $\delta \tilde{m}$ and $\delta \tilde{\lambda}$ are the solutions of 
\begin{eqnarray}
\label{eq:BLIncState}
\partial_t \delta \tilde{m} (t) + \widetilde{\nabla} \delta \tilde{m}(t) \star \tilde{v}_t + \widetilde{\nabla} \tilde{m}(t) \star \delta \tilde{v}(t) = 0\\
\label{eq:BLIncAdjoint}
-\partial_t \delta \tilde{\lambda}(t) - \widetilde{\nabla} \cdot ( \delta \tilde{\lambda}(t) \star \tilde{v}_t ) = 0
\end{eqnarray}
\noindent in the BL domain. The initial conditions are, respectively, $m(0) = I_0$, $\lambda(1) = -\frac{2}{\sigma^2}(m(1)-I_1)$, $\delta \tilde{m}(0) = 0$, 
$\delta \tilde{\lambda}(1) = -\frac{2}{\sigma^2} \delta \tilde{m}(1)$.

\subsubsection{BL PDE-constrained LDDMM based on the state equation}

This method departs from the original BL PDE-constrained LDDMM by using $m(t) = I_0 \circ \phi(t)$ and $\lambda(t) = J(t) \lambda(1) \circ \psi(t)$, 
where $\psi$ is the inverse of $\phi$, and $J$ is the Jacobian determinant of $\psi$~\cite{Hart_09,Hernandez_19c}.
The transformations $\phi$ and $\psi$ and the scalar field $J$ are computed from the inclusion of the truncated displacement fields ($\tilde{u}(t)$ and $\tilde{\nu}(t)$) 
and the corresponding Jacobian
\begin{eqnarray}
 \phi(t) = id - \iota(\tilde{u}(t)), \\ 
 \psi(t) = id - \iota(\tilde{\nu}(t)), \\ 
 J(t) = 1 - \iota(\tilde{U}(t)),
\end{eqnarray}
\noindent where
\begin{eqnarray}
\label{eq:BLDefState}
\partial_t \tilde{u}(t) + \tilde{D} \tilde{u}(t) \star \tilde{v}_t = \tilde{v}_t \\
-\partial_t \tilde{\nu}(t) - \tilde{D} \tilde{\nu}(t) \star \tilde{v}_t = -\tilde{\nu}_t \\
-\partial_t \tilde{U}(t) - \tilde{v}_t \star \widetilde{\nabla} \tilde{U}(t) = -\widetilde{\nabla \cdot} \tilde{v} + \tilde{U}(t) \star \widetilde{\nabla \cdot} \tilde{v}(t)
\end{eqnarray}

\noindent with initial conditions $\tilde{u}(0)= 0$, $\tilde{\nu}(1) = 0$, and $\tilde{U}(1) = 0$.

% \color{red}
% ... Cuidado que en SIAM creo que me deje por poner las ecuaciones incrementales ...
% \color{black}

The incremental state and adjoint variables are computed from
\begin{eqnarray}
\delta m(t) = \nabla I_0 \circ \phi(t) \cdot \iota(\delta \tilde{u}(t))\\
\delta \lambda(t) = J(t) \nabla \lambda(1) \circ \psi(t) \cdot \iota(\delta \tilde{\nu}(t)) % En el codigo esta con -
\end{eqnarray}
where $\delta \tilde{u}$ and $\delta \tilde{\nu}$ are computed from
\begin{eqnarray}
\label{eq:BLIncDefState}
\partial_t \delta \tilde{u}(t) + \tilde{D} \delta \tilde{u}(t) \star \tilde{v}_t + \tilde{D} \tilde{u}(t) \star \delta \tilde{v}_t = \delta \tilde{v}_t \\
-\partial_t \delta \tilde{\nu}(t) - \tilde{D} \delta \tilde{\nu}(t) \star \tilde{v}_t - \tilde{D} \tilde{\nu}(t) \star \delta \tilde{v}_t = -\delta \tilde{\nu}_t
\end{eqnarray}

\noindent with initial conditions $\delta \tilde{u}(0) = 0$ and $\delta \tilde{\nu}(1) = 0$.

\subsubsection{BL PDE-constrained LDDMM based on the deformation state equation}

The third method is formulated from the minimization of Equation~\ref{eq:BLEnergy}
subject to the truncated displacement state equation~\cite{Hernandez_19c}
\begin{equation}
\partial_t \tilde{u}(t) + \tilde{D} \tilde{u}(t) \star \tilde{v}_t = \tilde{v}_t.
\end{equation}

The gradient and the Hessian-vector product are given by the equations
\begin{eqnarray}
\widetilde{(\nabla_{\tilde{v}} E_{aug}(\tilde{v}))_t} = \tilde{\Momentum} \tilde{v}_t + \tilde{\rho}(t) - \widetilde{D} \tilde{u}(t) \star \tilde{\rho}(t) \\
\widetilde{ (H_{\tilde{v}} E_{aug}(\tilde{v}))_t } \delta \tilde{v} (t) = \tilde{\Momentum} \delta \tilde{v}(t) + \delta \tilde{\rho}(t) - \widetilde{D} \delta \tilde{u}(t) \star \tilde{\rho}(t), 
\end{eqnarray}

\noindent where the displacement state variable $\tilde{u}$, the adjoint variable $\tilde{\rho}$, and their incremental counterparts $\delta \tilde{u}$ and $\delta \tilde{\rho}$ are computed from
\begin{eqnarray}
\partial_t \tilde{u}(t) + \widetilde{D} \tilde{u}(t) \star \tilde{v}_t = \tilde{v}_t \\
\label{eq:BLDefAdjoint}
-\partial_t \tilde{\rho}(t) - \widetilde{\nabla \cdot}  (\tilde{\rho}(t) \star \tilde{v}_t ) = 0  \\
\partial_t \delta \tilde{u}(t) + \widetilde{D} \delta \tilde{u}(t) \star \tilde{v}_t + \widetilde{D} \tilde{u}(t) \star \delta \tilde{v}(t) = \delta \tilde{v}(t) \\
\label{eq:BLIncDefAdjoint}
-\partial_t \delta \tilde{\rho}(t) - \widetilde{\nabla \cdot} ( \delta \tilde{\rho}(t) \star \tilde{v}_t ) = 0 
\end{eqnarray}
\noindent with initial conditions $\tilde{u}(0) = 0$, $\tilde{\rho}(1)= \pi( -\frac{2}{\sigma^2} (m(1)-I_1) \nabla m(1))$, $\delta \tilde{u}(0) = 0$, 
and $\delta \tilde{\rho}(1) = \pi( -\frac{2}{\sigma^2} \delta m(1) \nabla m(1))$.

\section{Semi-Lagrangian Runge-Kutta integration}
\label{sec:SL}

\subsection{Semi-Lagrangian integration in a spatial domain}

Semi-Lagrangian (SL) integration methods~\cite{Staniforth_91} allow solving transport equations of the general form 

\begin{equation}
\label{eq:transport}
 D_t u = f(u,v),
\end{equation}

\noindent where $u: \Omega^d \times [0,1] \rightarrow \IR$ is a scalar or a vector function varying in time, 
and $$D_t u = \partial_t u + D u \cdot v.$$ 
SL methods combine the most interesting properties of Eulerian and Lagrangian schemes.
On the one hand, SL methods involve following the characteristic lines of the differential equation,
similarly to Lagrangian approaches.
On the other hand, the equation is solved on the regular grid, similarly to Eulerian approaches.
As a result, SL methods are unconditionally stable as Lagrangian schemes.
This means that the time sampling can be selected according to accuracy considerations rather than
stability considerations.
SL methods allow selecting a time sampling usually much smaller than Eulerian methods yielding a sensible 
reduction of the computational complexity.

SL schemes involve two steps. 
First, the departure points are computed solving the characteristic equation 

\begin{equation}
 D_t X(t) = v(t, X(t)),
\end{equation}

\noindent with initial condition $X(0) = x$.
The direction of the time integration can be forward or backward, 
depending on the direction of the time integration of the transport equation.
From the several methods proposed in the literature for solving the characteristic equation,
we use the approach given by Mang et al. in~\cite{Mang_17b}
\begin{eqnarray}
 X_* = x - \delta t \cdot v \\
 v_* = v \circ X_* \\
 X_* = x - 0.5 \textnormal{ } \delta t \cdot (v_* + v).
\end{eqnarray}

Second, the transport equation (Equation~\ref{eq:transport}) is solved in the Eulerian grid
\begin{equation}
\label{eq:TransportOmega}
 D_t u(X(t), t) = f(u(X(t),t), v(t, X(t))) 
\end{equation}

\noindent along the characteristic line $X$. 
The use of Runge-Kutta (RK) integration has been recently proposed in this step, yielding a higher-order accurate SL-RK method~\cite{Guo_13}. 
The velocity field needs to be estimated at points that do not belong to the Eulerian grid. 
Therefore, an interpolator is needed. Cubic interpolation is the method of choice
for SL schemes~\cite{Riishojgaard_98}.

\subsection{Semi-Lagrangian integration in a band-limited domain}

In $\tilde{\Omega}$, the transport equations are of the general form
\begin{equation}
\label{eq:BLtransport}
 \tilde{D}_t \tilde{u} = \tilde{f}(\tilde{u},\tilde{v}),
\end{equation}
\noindent where 
\begin{equation}
\tilde{D}_t \tilde{u} = \partial_t \tilde{u} + \tilde{D} \tilde{u} \star \tilde{v}.
\end{equation}

\noindent The characteristics are computed from 
\begin{equation}
 D_t X(t) = \iota(\tilde{v})(t, X(t)),
\end{equation}
\noindent and the transport equation is solved from
\begin{equation}
\label{eq:TransportOmegaTilde}
 \tilde{D}_t \tilde{u}(\tilde{X}(t), t) = \tilde{f}(\tilde{u}(\tilde{X}(t),t), \tilde{v}(t, \tilde{X}(t))). 
\end{equation}

\subsection{Semi-Lagrangian Runge-Kutta integration in PDE-LDDMM}

In this work, SL-RK integration has been implemented in $\Omega$ and $\tilde{\Omega}$ for the spatial and band-limited versions of the three 
PDE-LDDMM variants. To be able to apply SL integration, the differential equations need to be written in the shape of equations
~\ref{eq:TransportOmega} or~\ref{eq:TransportOmegaTilde}. We focus on the derivation for the BL domain $\widetilde{\Omega}$.
The derivation for the spatial domain can be performed analogously.

The state equations, the deformation state equations, and their incremental counterparts (Equations \ref{eq:BLState}, \ref{eq:BLDefState}, 
\ref{eq:BLIncState}, \ref{eq:BLIncDefState}) are already in the shape of equation~\ref{eq:TransportOmegaTilde} by just moving to
the right-hand-side of the equation a remaining term. For the adjoint and the incremental adjoint equations (Equations \ref{eq:BLAdjoint}, 
\ref{eq:BLDefAdjoint}, \ref{eq:BLIncAdjoint}, \ref{eq:BLIncDefAdjoint}) we use the identity 
\begin{eqnarray}
% \nabla \cdot u v = u \nabla \cdot v + v \cdot \nabla u \textnormal{ in } \Omega \\
\widetilde{\nabla \cdot} (\tilde{u} \star \tilde{v}) = \tilde{u} \widetilde{\nabla \cdot} \tilde{v} + \tilde{v} \star \widetilde{\nabla} \tilde{u} \textnormal{ in } \tilde{\Omega} 
\end{eqnarray}
\noindent and move the divergence term to the right-hand-side of the transformed equation. 
Table~\ref{table:DtEquations} gathers the expressions of the resulting differential equations, needed for the implementation of BL PDE-constrained LDDMM
methods in SL form.
For SL-RK, the right-hand side expressions can be directly plugged into an RK differential solver.

\setlength{\tabcolsep}{10pt} % Default value: 6pt
\renewcommand{\arraystretch}{1.5} % Default value: 1

\begin{table*}[!t]
\scriptsize
\centering
\begin{tabular}{|c|c|}
\hline
Equation in original form & Equation in SL form\\
\hline
 $\partial_t m(t) + \nabla m(t) \cdot \iota(\tilde{v}_t) = 0$ (eq.~\ref{eq:BLState}) & $D_t m(t) = 0$ \\
 $-\partial_t \lambda(t) - \nabla \cdot ( \lambda(t) \cdot \iota(\tilde{v}_t) ) = 0$ (eq.~\ref{eq:BLAdjoint}) & $-D_t \lambda(t) = \lambda(t) \nabla \cdot \iota(\tilde{v_t})$ \\
 $\partial_t \tilde{u}(t) + \widetilde{D} \tilde{u}(t) \star \tilde{v}_t = \tilde{v}_t$ (eq.~\ref{eq:BLDefState}) & $\tilde{D}_t \tilde{u}(t) = \tilde{v_t}$ \\
 $ -\partial_t \tilde{\rho}(t) - \widetilde{\nabla \cdot}  (\tilde{\rho}(t) \star \tilde{v}_t ) = 0 $ (eq.~\ref{eq:BLDefAdjoint}) & $-\tilde{D}_t \tilde{\rho}(t) = \tilde{\rho}(t) \star \widetilde{\nabla \cdot} \tilde{v_t}$ \\
 $\partial_t \delta \tilde{m} (t) + \widetilde{\nabla} \delta \tilde{m}(t) \star \tilde{v}_t + \widetilde{\nabla} \tilde{m}(t) \star \delta \tilde{v}(t) = 0$ (eq.~\ref{eq:BLIncState}) & $\tilde{D}_t \delta \tilde{m}(t) = -\widetilde{\nabla} \tilde{m}(t) \star \delta \tilde{v_t}$ \\
 $-\partial_t \delta \tilde{\lambda}(t) - \widetilde{\nabla} \cdot ( \delta \tilde{\lambda}(t) \star \tilde{v}_t ) = 0$ (eq.~\ref{eq:BLIncAdjoint}) & $-D_t \delta \tilde{\lambda}(t) = \tilde{\lambda}(t) \star \widetilde{\nabla \cdot} \tilde{v_t}$ \\
 $\partial_t \delta \tilde{u}(t) + \widetilde{D} \delta \tilde{u}(t) \star \tilde{v}_t + \widetilde{D} \tilde{u}(t) \star \delta \tilde{v}(t) = \delta \tilde{v}(t)$ (eq.~\ref{eq:BLIncDefState}) &  $\tilde{D}_t \delta \tilde{u}(t) = \delta \tilde{v}_t - \tilde{D}\tilde{u}(t) \star \delta \tilde{v}(t)$ \\
 $-\partial_t \delta \tilde{\rho}(t) - \widetilde{\nabla \cdot} ( \delta \tilde{\rho}(t) \star \tilde{v}_t ) = 0$ (eq.~\ref{eq:BLIncDefAdjoint}) & $-\tilde{D}_t \delta \tilde{\rho}(t) = \delta \tilde{\rho}(t) \star \widetilde{\nabla \cdot} \tilde{v_t}$ \\ 
\hline
\end{tabular}
\caption{Original PDEs involved in BL PDE-constrained LDDMM and corresponding PDEs written in SL form.}
\label{table:DtEquations}
\end{table*}

\section{Experimental setup}
\label{sec:ExpSetup}

In this work, we evaluate the performance of SL-RK integration in all the variants of PDE-constrained LDDMM (see Table~\ref{table:Methods}). 
The evaluation has been performed consistently with our previous work~\cite{Hernandez_19,Hernandez_19c}, in order to show the improvement of the 
proposed integration method over RK integration.
In addition, we have performed an extensive evaluation of the most memory efficient stationary methods 
in the frameworks of Klein et al.~\cite{Klein_09} and Rohlfing et al.~\cite{Rohlfing_12}
in order to establish the position achieved by PDE-constrained LDDMM methods in these evaluation rankings.

\begin{table*}[!t]
\begin{center}
\scriptsize
\begin{tabular}{|c|c|c|c|}
\hline
Abbreviation & Method & Integration & Publication \\
\hline
PDE-LDDMM RK & PDE-constrained LDDMM & RK & Mang et al. \cite{Mang_15} \\
\hline
PDE-LDDMM st. eq. RK & \begin{tabular}{c} PDE-constrained LDDMM \\ based on the state equation \end{tabular} & RK & Hernandez \cite{Hernandez_19c} \\   
\hline
PDE-LDDMM def. eq. RK & \begin{tabular}{c} PDE-constrained LDDMM \\ based on the deformation state equation \end{tabular} & RK & Hernandez \cite{Hernandez_19c} \\   
\hline
\hline
PDE-LDDMM SL & PDE-constrained LDDMM & SL-RK & Mang et al. \cite{Mang_17b} \\ 
\hline
PDE-LDDMM st. eq. SL & \begin{tabular}{c} PDE-constrained LDDMM \\ based on the state equation \end{tabular} & SL-RK & this work \\   
\hline
PDE-LDDMM def. eq. SL & \begin{tabular}{c} PDE-constrained LDDMM \\ based on the deformation state equation \end{tabular} & SL-RK & this work \\   
\hline
\hline
BL PDE-LDDMM RK & band-limited PDE-constrained LDDMM & RK & Hernandez \cite{Hernandez_19} \\ 
\hline
BL PDE-LDDMM st. eq. RK & \begin{tabular}{c}band-limited PDE-constrained LDDMM \\ based on the state equation \end{tabular} & RK & Hernandez \cite{Hernandez_19c} \\
\hline
BL PDE-LDDMM def. eq. RK & \begin{tabular}{c} band-limited PDE-constrained LDDMM \\ based on the deformation state equation \end{tabular} & RK & Hernandez \cite{Hernandez_19c} \\
\hline
\hline
BL PDE-LDDMM SL & band-limited PDE-constrained LDDMM & SL-RK & this work \\ 
\hline
BL PDE-LDDMM st. eq. SL & \begin{tabular}{c} band-limited PDE-constrained LDDMM \\ based on the state equation \end{tabular} & SL-RK &  this work \\
\hline
BL PDE-LDDMM def. eq. SL &  \begin{tabular}{c} band-limited PDE-constrained LDDMM \\ based on the deformation state equation \end{tabular} & SL-RK & this work \\
\hline
\end{tabular}
\caption{List of the PDE-constrained LDDMM methods compared in this work. The table gathers the abbreviations used in this experimental section, 
the name of the method given in Section 2, the integration scheme, and the publication where the method first appeared. }
\label{table:Methods}
\end{center}
\end{table*}

\subsection{Datasets}

\noindent We have used five different databases in our evaluation:

{\bf NIREP16} contains 16 skull-stripped brain images with the segmentation of 32 gray matter structures.
The dimension of the images is 256 $\times$ 300 $\times$ 256 with a voxel size of $ 0.7 \times 0.7 \times 0.7$ mm. 
The acquisition and post-processing details can be found at the web page~\url{http://www.nirep.org}.
The most remarkable features of this dataset are the excellent image quality and the ventricle sizes that are usually small.
The geometry of the segmentations provides a specially challenging framework for deformable registration evaluation.
% Calidad de las imagenes muy buena
% Ventriculos realmente pequeños
% Las segmentaciones van super por las circunvoluciones, curvatura

{\bf LPBA40} contains 40 skull-stripped brain images without the cerebellum and the brain stem.
LPBA40 is provided with the segmentation of 50 gray matter structures together with the caudate, putamen, and hippocampus.
% Image dimension is 180 $\times$ 216 $\times$ 180.
LPBA40 protocols can be found at

\noindent \url{http://www.loni.ucla.edu/Protocols/LPBA40}.
The image quality in LPBA40 is, overall, acceptable.
The variability of the ventricle sizes is high. 
% Calidad de las imagenes no tan buena como NIREP
% Uno en concreto malilla
% Mucha variabilidad de tamaño en los ventriculos, de los que no se da segmentacion
% Las segmentaciones son regiones a piñon, sin curvatura

{\bf IBSR18} contains 18 brain images with the segmentation of 96 cerebral structures.
The masks for skull-stripping are available with the dataset. 
In addition, the release {\bf IBSR\_V2.0} skull-stripped NIFTI~\cite{Rohlfing_12} contains 18 
skull-stripped brain images with the segmentation of 62 cerebral structures.
% Image dimension is 256 $\times$ 128 $\times$ 256, although in this work we cropped and 
% padded the images to size 220 $\times$ 180 $\times$ 220.
% 220x180x220 para Rohlfing
% 180x168x180 para Klein
This dataset provides the segmentation of brain structures of interest for the
evaluation of image registration methods.
The image quality is low. For example, most of the images show motion artifacts.
The variability of the ventricle sizes is high. 
% \color{red}
% In our work, we selected 12 out of the 18 images.
% The images were cropped and padded to size 220 $\times$ 180 $\times$ 220. 
% From the 62 labels, we selected 19 labels with the most important subcortical structures.
% ... Analizar los resultados. Conviene quitar esos pacientes si me quedo con las 19 etiquetas?
% Puede levantar ampollas en los revisores.
% \color{black}

{\bf CUMC12} contains 12 full brain images with the segmentation of 130 cerebral structures. % 132 menos 2 nan
The masks for skull-stripping are available with the dataset. 
% Image dimension is 256 $\times$ 124 $\times$ 256.
Overall, the image quality is acceptable, although some of the images are noisy.
The contrast of the images is low.
The variability of the ventricle sizes is high. 

{\bf MGH10} contains 10 full brain images with the segmentation of 106 cerebral structures. % 152 menos *** nan
The masks for skull-stripping are available with the dataset. 
% Image dimension is 182 $\times$ 218 $\times$ 182.
Overall, the image quality is acceptable, although some of the images are noisy.
The contrast of the images is low.
Ventricle sizes are usually all big. 

\subsection{Image registration pipeline}

The evaluation consistent with our previous work was performed in a subsampled NIREP16 database.
The registrations were carried out from the first subject to every other subject in the database, 
yielding to a total of 15 registrations per method.
The subsampled NIREP16 database was obtained from the resampling of the images into volumes of size 
$180 \times 210 \times 180$ with a voxel size of $1.0 \times 1.0 \times 1.0$ mm after the alignment 
to a common coordinate system using affine transformations. 
The images were scaled between $0$ and $1$.
The affine alignment and subsampling were performed using the Insight Toolkit (ITK).
The PDE-constrained registration methods were executed directly on this dataset.
% We run both single- and multi-resolution versions of the registration methods.
For benchmarking, we run single- and multi-resolution versions of the SyN version of ANTS diffeomorphic registration~\cite{Avants_08} 
with $L^2$ image similarity (ANTS-SSD).

The evaluation in the framework of Klein et al. was performed in NIREP16, LPBA40, IBSR18, CUMC12, and MGH10 databases.
The IBSR18, CUMC12 and MGH10 images normalized with respect to the MNI152 space were used as input data.
The registrations were carried out from every subject to every other subject in each database yielding to a total of 2154 registrations per method.
The evaluation in the framework of Rohlfing et al. was performed in IBSR18 database, with a total of 306 registrations per method.
% This removed any bias in the selection of source and target images in the registration evaluation.
The NIREP16, LPBA40, IBSR18, CUMC12, and MGH10 images were preprocessed similarly to~\cite{Klein_09}. 
In the first place, N4 bias field correction and histogram matching were applied to all the images.
To perform these preprocessing steps we used the algorithms available in ITK.
The images were scaled between $0$ and $1$.
Next, we performed an affine registration between all the image pairs.
Instead of using the affine registered images as input of our non-rigid registration methods, 
we used the affine transformation as input, and it was included in the parameterization of the
diffeomorphic transformations.
% The PDE-constrained registration methods were embedded into a multi-resolution scheme.
% The images were subsampled, and the velocity fields were resampled similarly to~\cite{Hernandez_17,Hernandez_18}.
% The PDE-constrained registration methods were executed on each resolution level.

Subsampled NIREP experiments were run on a cluster equipped with one NVidia Titan RTX with 24 GBS of video memory
and an Intel Core i7 with 64 GBS of DDR3 RAM.
NIREP16, LPBA40, IBSR18, CUMC12, and MGH10 experiments were run on a cluster equipped with four NVidia GeForce GTX 1080 ti with 11 GBS of video 
memory and an Intel Core i7 with 64 GBS of DDR3 RAM.
The codes were developed in the GPU with Matlab 2017a and Cuda $8.0$.
Since Matlab lacks a 3D GPU cubic interpolator, we implemented in a Cuda MEX file the GPU cubic interpolator 
with prefiltering proposed in~\cite{Ruijters_12}.

\subsection{Parameter configuration}

Regularization parameters were selected from a search of the optimal parameters in the registration experiments 
performed in our previous work~\cite{Hernandez_19c}.
We selected the parameters $\sigma^2 = 1.0$, $\alpha = 0.0025$, and $s = 2$. 

% For the single-resolution experiments, the optimization was run a maximum of 10 iterations with the stopping 
% conditions used in~\cite{Mang_15}.
The optimization was run a maximum of 10 iterations with the stopping conditions used in~\cite{Mang_15}.
The maximum number of PCG iterations was selected equal to 5.
These parameters were selected as optimal in our previous work since the methods achieved state of the art 
accuracy at a reasonable amount of time~\cite{Hernandez_19}.
% For the multi-resolution experiments, the pyramid was built with 3 levels with the same number
% of outer and inner iterations than for the single-resolution.  

The experiments were performed with band sizes of $32 \times 32 \times 32$ for BL PDE-constrained LDDMM based on 
the state and on the deformation state equations, and band sizes of $40 \times 40 \times 40$ for original BL PDE-constrained LDDMM.
This selection was found as optimal for each method in our previous work~\cite{Hernandez_19,Hernandez_19c}.

For SL-RK integration, $n_t$ was selected equal to $5$ for all the methods.
For RK integration, $n_t$ was selected equal to $25$ for the BL PDE-constrained LDDMM based on the state and on the 
deformation state equations, and $50$ for the spatial methods due to stability issues.
In the evaluation with LPBA40, IBSR18, CUMC12, and MGH10 datasets, $n_t=25$ showed stability issues in a considerable 
number of experiments and it was raised to $50$.

ANTS-SSD was run with the following parameters for the single-resolution experiments 

\noindent \texttt{\$synconvergence="[50,1e-6,10]"},

\noindent \texttt{\$synshrinkfactors="1"}, 

\noindent and \texttt{\$synsmoothingsigmas="3vox"}.

For the multi-resolution experiments the parameters were set to

\noindent \texttt{\$synconvergence="[50x50x50,1e-6,10]"}, 

\noindent \texttt{\$synshrinkfactors="4x2x1"}, 

\noindent and \texttt{\$synsmoothingsigmas="3x2x1vox"}.

\noindent The selection of the number of iterations was in agreement with the number of outer $\times$ inner iterations used in Gauss--Newton--Krylov 
optimization.

\section{Results}
\label{sec:Results}

\subsection{Subsampled NIREP16 evaluation results}

\subsubsection{Convergence analysis}

% -- Sauron --

\begin{table*}[!t]
\begin{center}
\scriptsize
% \begin{tabular}{|p{55 mm}|c|c|c|c|c|}
% Single-resolution
% \\
\begin{tabular}{|c|c|c|c|c|c|}
\hline
Method & $n_t$ & $MSE_{rel}$ & $\Vert g\Vert_{\infty,{rel}}$ & $\min(J(\phi_{1}^v)^{-1})$ & $\max(J(\phi_{1}^v)^{-1})$ \\
\hline
St. PDE-LDDMM RK & 50 & 18.29 $\pm$ 2.83 & 0.07 $\pm$ 0.05 & 0.16 $\pm$ 0.05 & 3.70 $\pm$ 0.51 \\   
St. PDE-LDDMM st. eq. RK & 30 & 18.42 $\pm$ 2.71 & 0.18 $\pm$ 0.19 & 0.07 $\pm$ 0.06 & 2.95 $\pm$ 0.41 \\
St. PDE-LDDMM def. st. eq. RK & 30 & 16.83 $\pm$ 1.53 & 0.16 $\pm$ 0.08 & 0.12 $\pm$ 0.05 & 4.12 $\pm$ 1.41 \\
\hline
St. PDE-LDDMM SL & 5 & 19.55 $\pm$ 1.76 & 0.14 $\pm$ 0.13 & 0.03 $\pm$ 0.04 & 3.25 $\pm$ 0.29 \\ % Jac RK
St. PDE-LDDMM st. eq. SL & 5 & 18.34 $\pm$ 1.68 & 0.23 $\pm$ 0.21 & 0.03 $\pm$ 0.04 & 3.46 $\pm$ 0.69 \\ % Jac RK
St. PDE-LDDMM def. st. eq. SL & 5 & 17.10 $\pm$ 1.50 & 0.12 $\pm$ 0.05 & 0.14 $\pm$ 0.05 & 5.02 $\pm$ 1.01 \\ % Jac RK
\hline
NSt. PDE-LDDMM RK & 50 & 22.68 $\pm$ 5.49 & 0.26 $\pm$ 0.12 & 0.16 $\pm$ 0.05 & 3.38 $\pm$ 0.90 \\
NSt. PDE-LDDMM, st. eq. RK & 25$^*$ & 29.92 $\pm$ 2.41 & 0.32 $\pm$ 0.10 & 0.11 $\pm$ 0.07 & 2.93 $\pm$ 0.78 \\
NSt. PDE-LDDMM, def. st. eq. RK & 30 & 16.10 $\pm$ 1.70 & 0.20 $\pm$ 0.08 & 0.09 $\pm$ 0.03 & 4.55 $\pm$ 0.69 \\
\hline
NSt. PDE-LDDMM SL & 5 & 22.44 $\pm$ 6.26 & 0.51 $\pm$ 0.18 & 0.02 $\pm$ 0.02 & 3.96 $\pm$ 0.65 \\ % Jac RK
NSt. PDE-LDDMM st. eq. SL & 5 & 21.82 $\pm$ 2.51 & 0.49 $\pm$ 0.19 & 0.01 $\pm$ 0.02 & 4.20 $\pm$ 0.94 \\ % Jac RK
NSt. PDE-LDDMM def. st. eq. SL & 5 & 16.86 $\pm$ 1.91 & 0.13 $\pm$ 0.04 & 0.09 $\pm$ 0.02 & 7.35 $\pm$ 2.63 \\ % Jac RK
\hline
BL St. PDE-LDDMM 40x RK & 25 & 20.99 $\pm$ 2.59 & 0.04 $\pm$ 0.03 & 0.21 $\pm$ 0.04 & 3.32 $\pm$ 0.41 \\
BL St. PDE-LDDMM st. eq. 32x RK & 25 & 18.53 $\pm$ 1.71 & 0.02 $\pm$ 0.01 & 0.11 $\pm$ 0.07 & 3.10 $\pm$ 0.35 \\
BL St. PDE-LDDMM def. st. eq. 32x RK & 25 & 17.32 $\pm$ 1.68 & 0.03 $\pm$ 0.01 & 0.13 $\pm$ 0.05 & 3.76 $\pm$ 0.49 \\ 
\hline
BL St. PDE-LDDMM 40x SL & 5 & 20.41 $\pm$ 1.89 & 0.02 $\pm$ 0.00 & 0.32 $\pm$ 0.02 & 9.65 $\pm$ 4.77  \\
BL St. PDE-LDDMM st. eq. 32x SL & 5 & 19.89 $\pm$ 1.76 & 0.01 $\pm$ 0.00 & 0.29 $\pm$ 0.03 & 7.94 $\pm$ 2.81 \\
BL St. PDE-LDDMM def. st. eq. 32x SL & 5 & 17.77 $\pm$ 1.66 & 0.04 $\pm$ 0.01 & 0.21 $\pm$ 0.03 & 8.23 $\pm$ 3.24 \\                               
\hline
BL NSt. PDE-LDDMM 40x RK & 25 & 29.30 $\pm$ 3.50 & 0.21 $\pm$ 0.08 & 0.24 $\pm$ 0.05 & 2.25 $\pm$ 0.51 \\
BL NSt. PDE-LDDMM, st. eq. 32x RK & 25 & 29.21 $\pm$ 3.96 & 0.24 $\pm$ 0.11 & 0.15 $\pm$ 0.06 & 2.71 $\pm$ 0.48 \\
BL NSt. PDE-LDDMM, def. st. eq. 32x RK & 25 & 15.68 $\pm$ 1.52 & 0.04 $\pm$ 0.01 & 0.09 $\pm$ 0.04 & 4.77 $\pm$ 0.78 \\
\hline
BL NSt. PDE-LDDMM 40x SL & 5 & 19.14 $\pm$ 2.06 & 0.11 $\pm$ 0.05 & 0.01 $\pm$ 0.07 & 5.25 $\pm$ 0.80 \\ 
BL NSt. PDE-LDDMM st. eq. 32x SL & 5 & 20.49 $\pm$ 1.75 & 0.17 $\pm$ 0.15 & 0.07 $\pm$ 0.03 & 4.51 $\pm$ 0.61 \\
BL NSt. PDE-LDDMM def. st. eq. 32x SL & 5 & 16.80 $\pm$ 1.57 & 0.07 $\pm$ 0.03 & 0.09 $\pm$ 0.03 & 7.37 $\pm$ 1.43 \\
\hline
\end{tabular}

\caption{ Subsampled NIREP16. Convergence results.
Mean and standard deviation of the relative image similarity error expressed in \% ($MSE_{rel}$), 
the relative gradient magnitude ($\Vert g\Vert_{\infty,{rel}}$), 
and the Jacobian determinant extrema associated with the transformation $(\phi^{v}_{1})^{-1}$.
St. stands for the stationary and NSt. for the non-stationary parameterization of diffeomorphisms.
(*) The $n_t$ selected for NSt. PDE-LDDMM, st. eq. with RK integration was due to a memory load greater than 
the GPU maximum capacity (24 GBS). }
\label{table:QuantitativeResults}
\end{center}
\end{table*}

%------------------------------------------------------------------------- 

Table~\ref{table:QuantitativeResults} shows, averaged by the number of experiments, 
the mean and standard deviation of the relative image similarity error after registration, 
$$MSE_{rel} = \frac{\Vert m(1) - I_1 \Vert_{L^2}^2}{\Vert I_0 - I_1 \Vert_{L^2}^2},$$
\noindent the relative gradient magnitude, 
$$\Vert g \Vert_{\infty,rel} = \frac{\Vert \nabla_{v} E({v}^n) \Vert_\infty}{\Vert \nabla_{v} E({v}^0) \Vert_\infty},$$ 
and the extrema of the Jacobian determinant obtained with PDE-constrained LDDMM in the subsampled NIREP16 dataset.

We have grouped the analysis of the $MSE_{rel}$ values by integration scheme, variant, image domain, and diffeomorphism
parameterization:

\begin{itemize}
 \item RK vs SL-RK. The $MSE_{rel}$ value achieved by SL-RK methods is close or improves RK methods.
 It drives our attention the bad performance of PDE-LDDMM st. eq. for the non-stationary parameterization and RK integration
 which is indeed improved by SL-RK integration.
 \item Variants. The best performing variant is PDE-LDDMM def. st. eq. This result is persistent for the different
 integration schemes, the spatial or BL parameterization, and the diffeomorphic parameterization.
 \item SP vs BL. The $MSE_{rel}$ values obtained with the spatial methods are slightly degraded by the BL methods as expected.
 The degradation is only shown in some cases.
 \item St. vs NSt. The non-stationary methods do not outperform the stationary methods in a consistent manner.
 The (out)performance depends on the variant.
\end{itemize}

For the bad-performing methods, the high values of the relative gradient indicate the stagnation of the convergence.
For the BL methods, the relative gradient was reduced to average values ranging from $0.02$ to $0.04$, which means that the optimization 
was stopped in acceptable energy values in all cases.
All the Jacobians remained above zero. 

\subsubsection{Evaluation}

The evaluation is based on the accuracy of the registration results for template-based segmentation.
% We use the manual segmentations provided with the NIREP database as the gold standard.
The Dice Similarity Coefficient (DSC) is selected as the evaluation metric.
Given two segmentations $S$ and $T$, the DSC is defined as
\begin{equation}
DSC(S,T) = \frac{2 \textnormal{Vol} (S \cap T)}{\textnormal{Vol}( S ) + \textnormal{Vol}( T )}.
\end{equation}
\noindent This metric provides the value of 1 if $S$ and $T$ exactly overlap and gradually decreases towards 0 
depending on the overlap of the two volumes.
% As a reference, the average DSC achieved by SyN diffeomorphic registration~\cite{Avants_08} with the subsampled NIREP16 dataset 
% is $0.56$ for the SSD image similarity.
% is $0.60$ for Cross-Correlation and $0.56$ for SSD image similarities.

\begin{figure*} [!t]
\centering
\scriptsize
\begin{tabular}{c}
\includegraphics[angle = 0, width = 1.25 \textwidth]{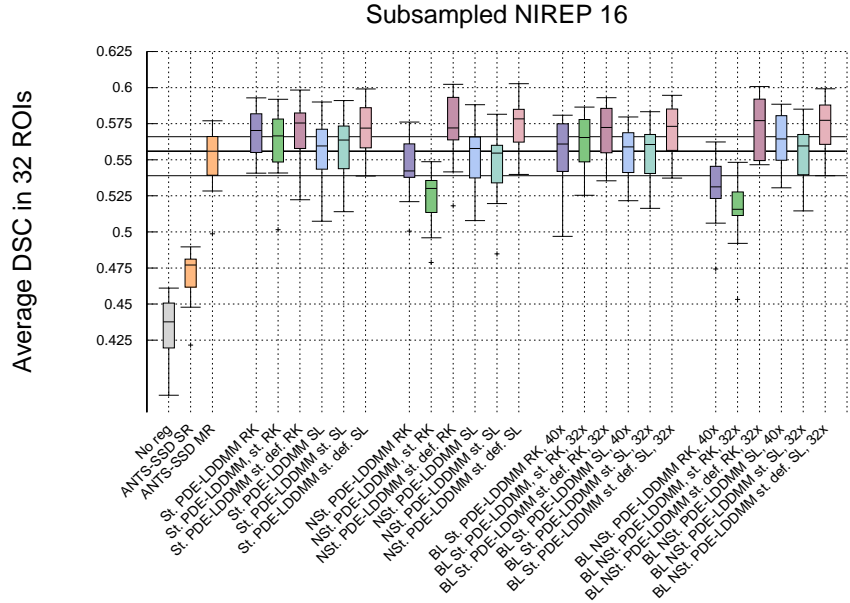}
\end{tabular}
\caption{Subsampled NIREP16. 
Volume overlap obtained by the registration methods measured in terms of the DSC between
the warped and the corresponding manual target segmentations. 
Box and whisker plots show the distribution of the DSC values averaged over the 32 NIREP manual segmentations.
The whiskers indicate the minimum and maximum of the DSC values.
The horizontal lines in the plot indicate the first, second, and third quartiles of multi-resolution ANTS-SSD.
} 
\label{fig:NIREPSubDSC}
\end{figure*}

Figure~\ref{fig:NIREPSubDSC} shows, in the shape of box and whisker plots, the statistical distribution of the 
DSC values obtained after the registration across the 32 segmented structures.
For the single-resolution experiments, the performance of the benchmark method ANTS-SSD was under $50 \%$.
For the multi-resolution experiments, the average DSC value achieved by ANTS-SSD equals to $55.59 \%$.
We have selected this value as a baseline of good registration accuracy for methods with $L^2$-based image 
similarity. 

A great number of PDE-LDDMM methods showed similar values or even improved ANTS-SSD performance.
The best performing variant was our PDE-LDDMM based on the deformation state equation (boxes in pink tones). 
This variant showed similar results for RK and SL-RK integration regardless the image domain and 
diffeomorphism parameterization.

For the variant associated with the PDE-constrained benchmark methods~\cite{Mang_15,Mang_17b} (boxes in blue tones), 
RK integration slightly outperformed SL-RK integration for the stationary parameterization.
For the non-stationary parameterization the median DSC value for RK integration was under the value 
for SL-RK integration.
Similarly, RK slightly outperformed SL-RK integration for BL PDE-LDDMM and BL PDE-LDDMM based on the state 
equation for the stationary parameterization.
On the contrary, SL-RK integration greatly outperformed RK integration for the non-stationary parameterization.

The variant based on the state equation (boxes in green tones) performed similarly to benchmark PDE-LDDMM variant for the 
stationary parameterization for the same integration scheme. 
However, it is remarkable the low performance achieved by this variant for the non-stationary parameterization
and RK integration in both image domains.

\subsubsection{Computational complexity}

\begin{table*}[!tpb]

\begin{center}
\scriptsize

\begin{tabular}{|c|c|c|c|}
\hline
Method & Integration & VRAM (MBS) & total time (s) \\
\hline
St. PDE-LDDMM & RK & 10311 & 2281.40 $\pm$ 378.56 \\   
St. PDE-LDDMM st. eq. & RK & 17229 & 877.85 $\pm$ 176.03 \\
St. PDE-LDDMM def. st. eq. & RK & 15321 & 1298.50 $\pm$ 151.19 \\
\hline
St. PDE-LDDMM & SL-RK & 3901 & 118.29 $\pm$ 4.77 \\
St. PDE-LDDMM st. eq. & SL-RK & 6267 & 139.91 $\pm$ 22.06 \\
St. PDE-LDDMM def. st. eq. & SL-RK & 5935 & 221.01 $\pm$ 2.40 \\
\hline
NSt. PDE-LDDMM & RK & 14913 & 3443.17 $\pm$ 932.45 \\
NSt. PDE-LDDMM, st. eq. & RK & 20635 & 816.79 $\pm$ 480.61 \\
NSt. PDE-LDDMM, def. st. eq. & RK & 19065 & 2222.04 $\pm$ 369.96 \\
\hline
NSt. PDE-LDDMM & SL-RK & 11155 & 191.96 $\pm$ 45.53 \\
NSt. PDE-LDDMM st. eq. & SL-RK & 12429 & 214.09 $\pm$ 48.54 \\
NSt. PDE-LDDMM def. st. eq. & SL-RK & 12065 & 355.47 $\pm$ 4.38 \\
\hline
\hline
BL St. PDE-LDDMM 40x & RK & 5709 & 312.45 $\pm$ 5.22 \\
BL St. PDE-LDDMM st. eq. 32x & RK & 4743 & 315.00 $\pm$ 2.73 \\
BL St. PDE-LDDMM def. st. eq. 32x & RK & 1819 & 377.56 $\pm$ 4.26 \\
\hline
BL St. PDE-LDDMM 40x & SL-RK & 2685 & 68.66 $\pm$ 1.32 \\
BL St. PDE-LDDMM st. eq. 32x & SL-RK & 2877 & 80.76 $\pm$ 0.35 \\
BL St. PDE-LDDMM def. st. eq. 32x & SL-RK & 2365 & 116.99 $\pm$ 1.04 \\
\hline
BL NSt. PDE-LDDMM 40x & RK & 6657 & 825.15 $\pm$ 19.93 \\
BL NSt. PDE-LDDMM, st. eq. 32x & RK & 4789 & 535.52 $\pm$ 279.65 \\
BL NSt. PDE-LDDMM, def. st. eq. 32x & RK & 1863 & 905.48 $\pm$ 25.13 \\
\hline
BL NSt. PDE-LDDMM 40x & SL-RK & 6131 & 158.79 $\pm$ 15.29 \\
BL NSt. PDE-LDDMM st. eq. 32x & SL-RK & 6171 & 159.49 $\pm$ 20.59 \\
BL NSt. PDE-LDDMM def. st. eq. 32x & SL-RK & 5577 & 249.60 $\pm$ 10.85 \\
\hline

\end{tabular}
\\
\vspace{0.1 cm}
\caption{Computational complexity. GPU peak memory usage and mean and standard deviation of the total computation 
time. Experiments run in an NVidia Titan RTX with 24 GBS of video memory.
% per iteration and overall
% It should be noticed that 1 iteration of a method based on Gauss-Newton-Krylov optimization is comparable
% with 5 iterations of a MCC benchmark method based on gradient-descent.
}
\label{table:Complexity}

\end{center}

\end{table*}

Table~\ref{table:Complexity} shows the VRAM peak memory reached through the computations, 
and the average and standard deviation of the total computation time in the subsampled NIREP16 experiments.
For the spatial methods, SL-RK integration achieved a substantial time and memory reduction over RK integration. % x31.18
The time and memory reduction achieved by SL-RK over RK integration was also considerable for the BL parameterized methods.  % x3.99 x4.87 x2.41
For the stationary parameterization, the BL parameterization further decreased the complexity of spatial SL-RK integration methods, 
as expected. % x1.31 x1.63 x1.57
However, SL-RK integration did not reduced the memory VRAM usage for the non-stationary parameterization while the total computation 
time was effectively reduced.

\subsubsection{Qualitative assessment}

For a qualitative assessment of the proposed registration methods, we show the registration results obtained by 
Mang et al. benchmark methods~\cite{Mang_15,Mang_17}, and 
PDE-LDDMM based on the deformation state equation in a selected experiment representative of a difficult deformable registration problem.
For the non-stationary parameterization the images for a qualitative assessment were similar to the stationary parameterization.
Figure~\ref{fig:QualitativeBenchmark} shows the warped images, the difference between the warped and the target images after registration, 
and the velocity field.

\subsection{NIREP16, LPBA40, IBSR18, CUMC12, and MGH10 evaluation results}

\begin{figure*} [!t]
\centering
\scriptsize
\begin{tabular}{c}
\includegraphics[angle = 0, width = 0.60 \textwidth]{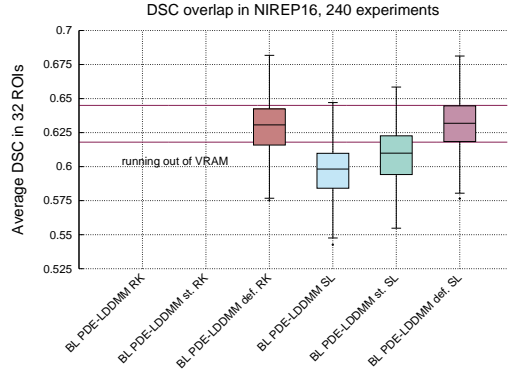}
\\
\end{tabular}
\caption{ NIREP16. Distribution of the DSC values averaged over the 32 NIREP manual segmentations in the 240 experiments.
The whiskers indicate the minimum and maximum of the DSC values.
The methods running out of 11 GBS VRAM are indicated in the plots.
The horizontal red lines indicate the first and third quartiles of BL PDE-LDDMM based on the deformation state equation.} 
\label{fig:NIREPDSC}
\end{figure*}

\begin{figure*}[t]
\centering
\scriptsize
\begin{tabular}{c}
 \includegraphics[angle = 0, width=0.75\textwidth]{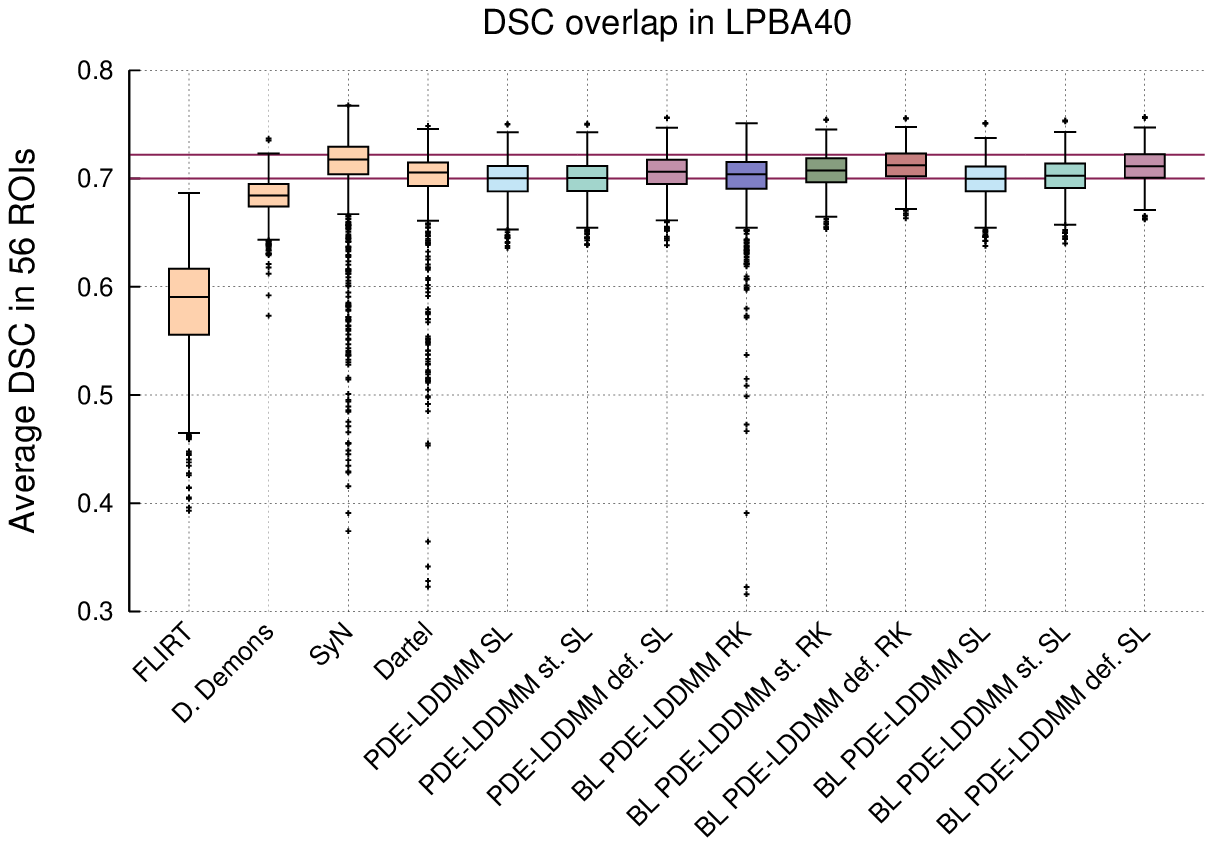} \\
 \includegraphics[angle = 0, width=0.75\textwidth]{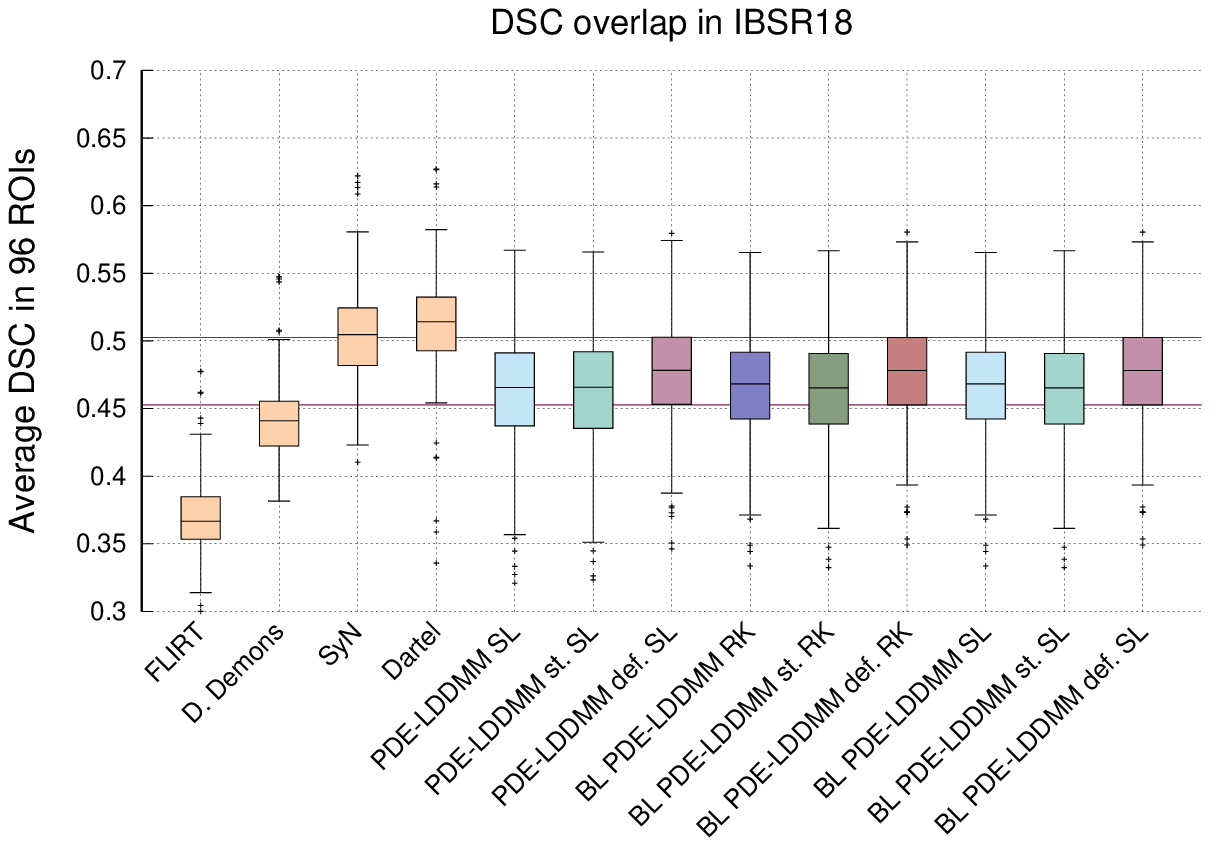} \\
\end{tabular}
\caption{ LPBA40 and IBSR18. 
Distribution of the DSC values averaged over the manual segmentations in the registration experiments.
The whiskers indicate the minimum and maximum of the DSC values.
The horizontal red lines indicate the first and third quartiles of BL PDE-LDDMM based on the deformation state equation.
} 
\label{fig:KleinDSC1}
\end{figure*}

\begin{figure*}[t]
\centering
\scriptsize
\begin{tabular}{c}
 \includegraphics[angle = 0, width=0.75\textwidth]{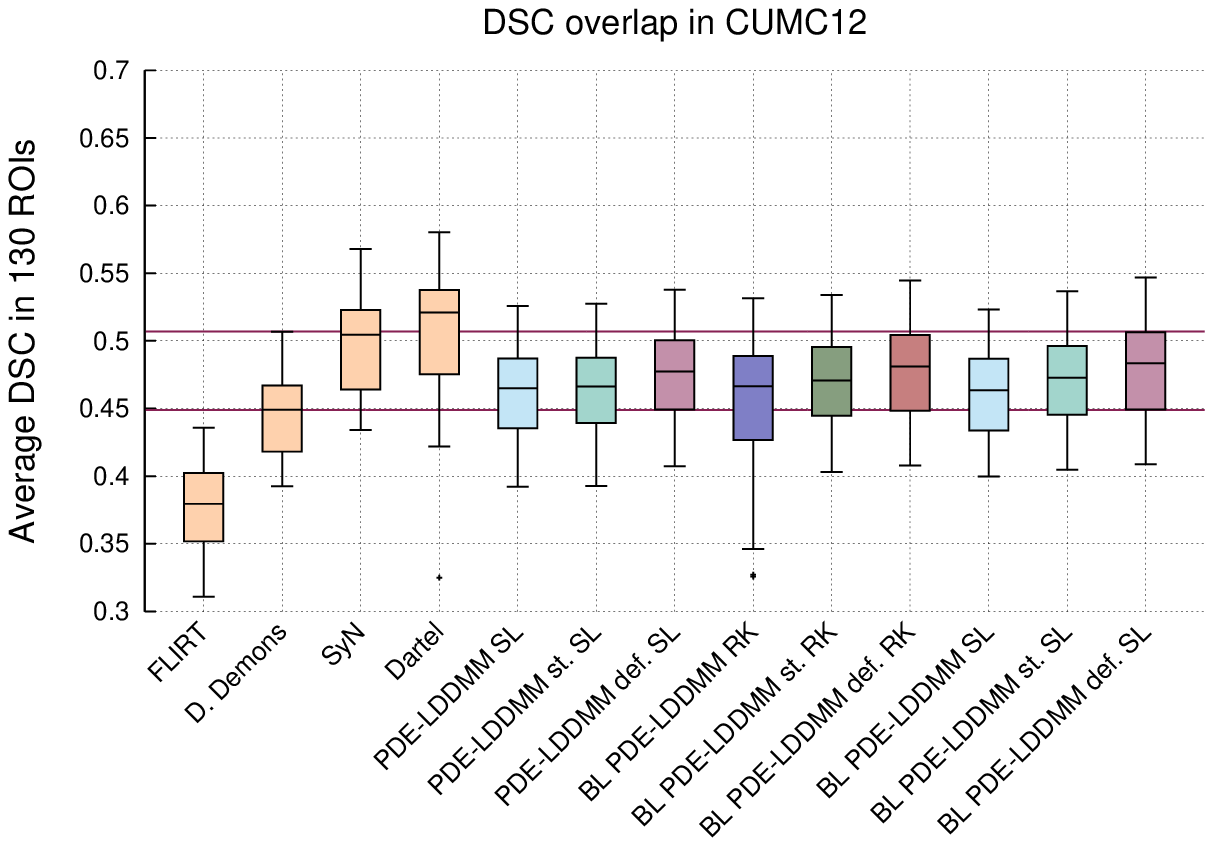} \\
 \includegraphics[angle = 0, width=0.75\textwidth]{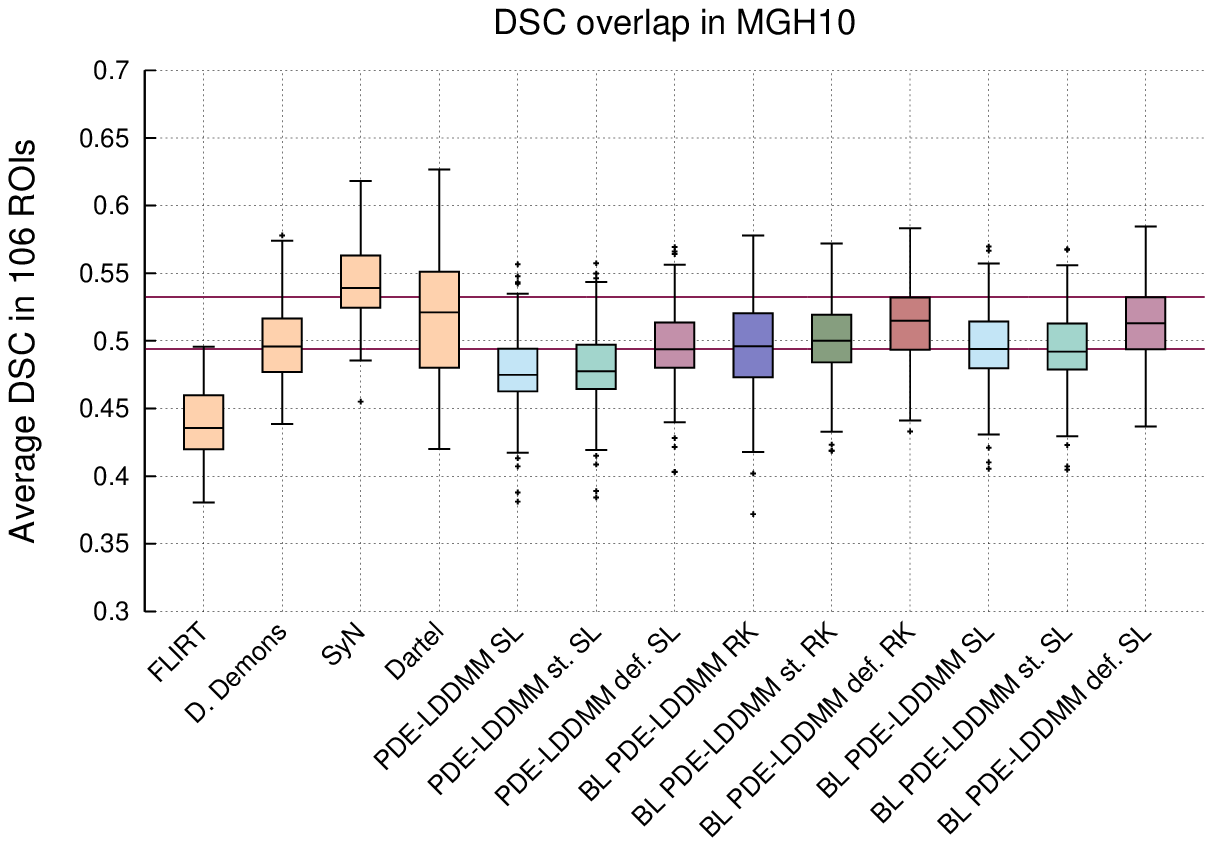} \\ 
\end{tabular}
\caption{ CUMC12 and MGH10. 
Distribution of the DSC values averaged over the manual segmentations in the registration experiments.
The whiskers indicate the minimum and maximum of the DSC values.
The horizontal red lines indicate the first and third quartiles of BL PDE-LDDMM based on the deformation state equation.
} 
\label{fig:KleinDSC2}
\end{figure*}

From the evaluation measurements used in Klein et al. framework, we focus on the accuracy of the registration results for 
template-based segmentation. 
Since~\cite{Klein_09}, this has been adopted as a widely extended criterion for non-rigid registration evaluation.
From the metrics proposed in~\cite{Klein_09}, we select the Dice Similarity Coefficient (DSC) as evaluation metric.
% Although other evaluation metrics have been used in the literature, the relative performance of the algorithms is little
% affected by the choice of the metric.
Figures~\ref{fig:NIREPDSC},~\ref{fig:KleinDSC1} and~\ref{fig:KleinDSC2} show the statistical distribution of the DSC values obtained after the registration 
across the manually segmented structures for the five databases.
For the NIREP16 dataset, we show the results obtained with ANTS-SSD. 
For the remaining databases, we include the results reported in~\cite{Klein_09} for affine registration (FLIRT), 
and three diffeomorphic registration methods: Diffeomorphic Demons, SyN, and Dartel.

On the one hand, the results from NIREP16 show how PDE-LDDMM based on the deformation state equation outperformed the other variants 
of PDE-constrained LDDMM. 
The distribution of the method parameterized in the BL domain resulted almost identical to the distribution of the method parameterized 
in the spatial domain.
As happened with subsampled NIREP16 evaluation, ANTS-SSD was among the worst performing methods.
BL PDE-LDDMM and BL PDE-LDDMM based on the state equation with RK integration exceeded the maximum VRAM capacity of our GPUs for this dataset.

On the other hand, the results from LPBA40, IBSR18, CUMC12, and MGH10 show that, from the PDE-constrained LDDMM methods, 
the best performing method was BL PDE-LDDMM based on the deformation state equation and SL-RK integration.
The performance of the method with RK integration was slightly lower.
The performance of original PDE-LDDMM, PDE-LDDMM based on the state equation and their BL versions was significantly lower in IBSR18, 
CUMC12, and MGH10 databases. 
For these methods, RK integration performed slightly better than SL-RK integration.
These results are consistent with NIREP16 evaluation results.

In IBSR18, CUMC12, and MGH10 databases, our PDE-constrained LDDMM methods were not able to reach SyN or Dartel performance.
This is probably because the image similarity metrics used in these methods (Cross-Correlation and multinomial model, respectively) 
favour the accuracy in template-based segmentation.
In contrast, PDE-constrained LDDMM uses SSD, which is known to restrict the performance in template-based segmentation.
However, in LPBA40 databases, our best performing PDE-constrained LDDMM methods overpass Dartel and almost reached SyN 
performance, while showing a significantly reduced number of outliers.

Our methods significantly outperformed FLIRT and Diffeomorphic Demons, where the third quartile in the distribution of our 
best performing method was close to the median of Demons for the four databases.
It should be noticed that Diffeomorphic Demons also uses SSD as image similarity metric.

\subsection{IBSR18 V2.0 evaluation results}

\begin{figure*}[h]
\centering
\scriptsize
\begin{tabular}{c}
 \includegraphics[angle = 0, width=0.75\textwidth]{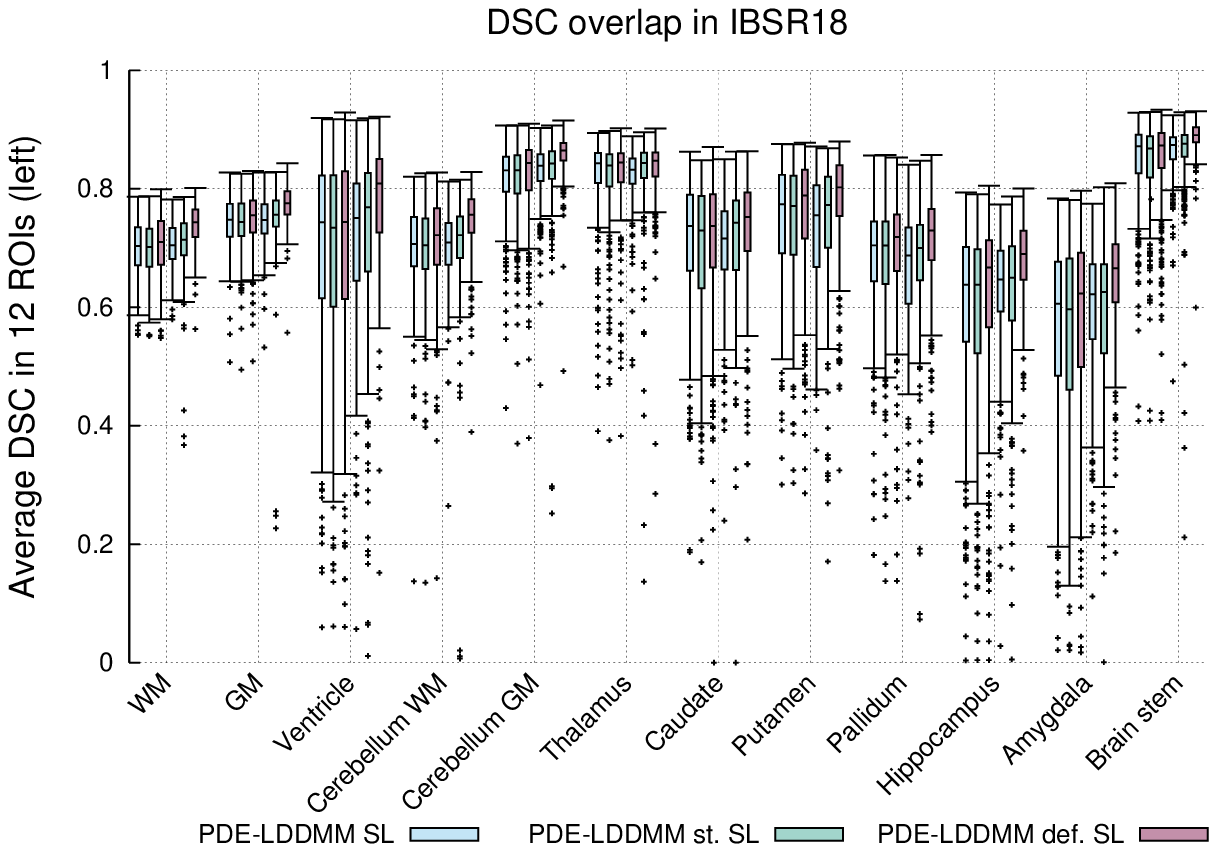} \\
 \includegraphics[angle = 0, width=0.75\textwidth]{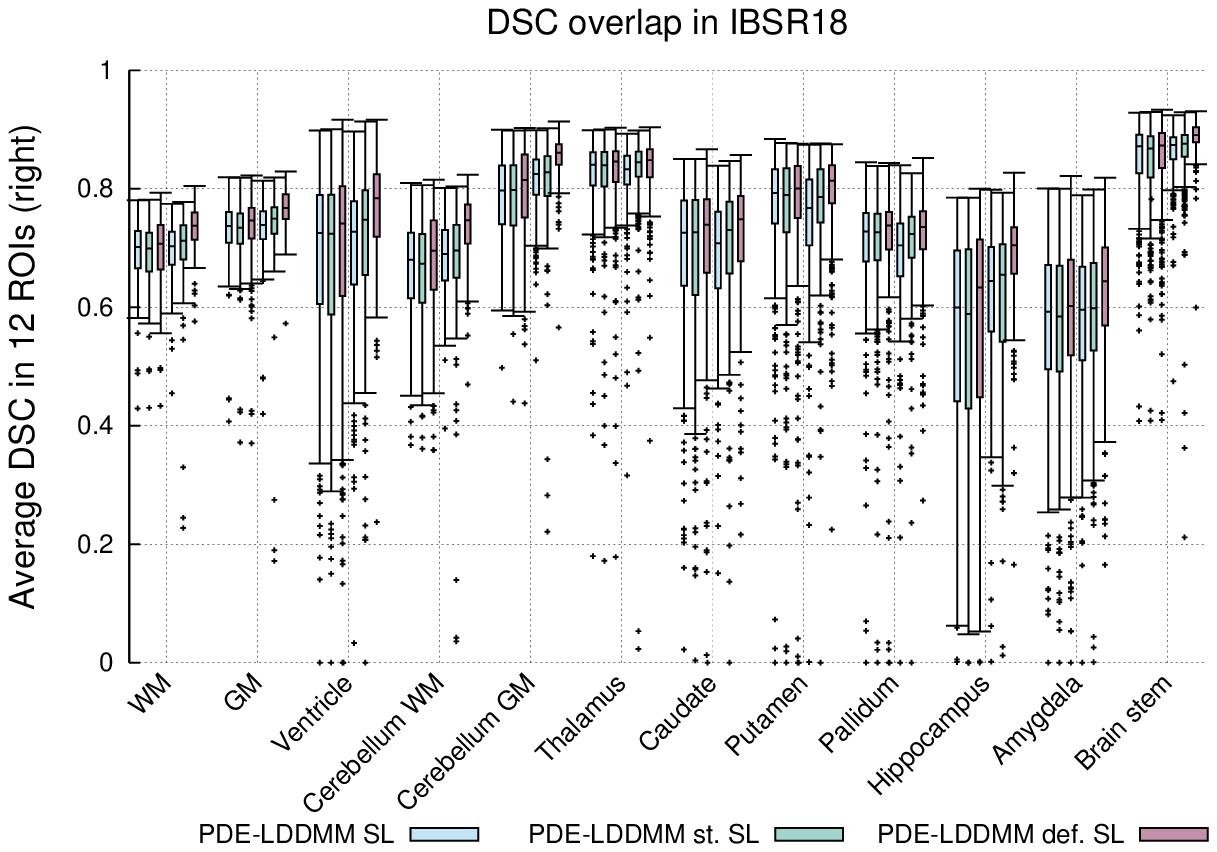} \\
\end{tabular}
\caption{ IBSR18 V2.0. Volume overlap obtained by the proposed registration methods.
The box and whisker plots show the distribution of the DSC values averaged over the manual segmentations for each region.
The whiskers indicate the minimum and maximum of the DSC values.
For each group of plots, the first three correspond with the spatial domain parameterization and the three last 
with the band-limited domain parameterization.
WM and GM stand for white matter and grey matter, respectively.
} 
\label{fig:IBSRDSC}
\end{figure*}

Figure~\ref{fig:IBSRDSC} shows the statistical distribution of the DSC values obtained by our proposed registration
methods in the regions of interest of Rohlfing et al. evaluation framework~\cite{Rohlfing_12}.
Consistently with the rest of the evaluation results, the best performing method was BL PDE-LDDMM based on the deformation 
state equation, which significantly outperformed the others in the great majority of regions.
% Figure~\ref{fig:Ventricles} shows the quality of the non-rigid alignment achieved at the ventricles by BL PDE-LDDMM based 
% on the deformation state equation in a representative experiment.

% \begin{figure*} [!t]
% \centering
% \small
% \begin{tabular}{cc}
% \includegraphics[angle = 0, width = 0.25 \textwidth]{./Figures/SourceTarget.eps}
% &
% \includegraphics[angle = 0, width = 0.25 \textwidth]{./Figures/PolzinTarget.eps} 
% \end{tabular}
% 
% \caption{ \small IBSR18 V2.0. Ventricle shapes before and after registration with BL PDE-LDDMM based on the deformation state
% equation with SL-RK integration.
% The red ventricles are generated from the manual segmentation of the target image.
% The green ventricles are generated from the manual segmentation of the source image.
% The pink shape is obtained warping the ventricles of the source image with the transformation estimated by our method.
% } 
% \label{fig:Ventricles}
% \end{figure*}

\section{Discussion}
\label{sec:Discussion}

The increase of the computational efficiency achieved by the combined over the split BL and SL-RK approaches was significant in terms of 
computation time and memory. 
The reduction of the memory requirements allowed us to perform the evaluation of the SL-RK PDE-constrained LDDMM methods extensively, 
even in the highest-resolution level of NIREP16.

In all the evaluation frameworks, BL PDE-constrained LDDMM based on the deformation state equation with SL-RK integration resulted our best performing method.
This method achieved an identical DSC distribution compared with RK integration in the NIREP16 database.
In NIREP16 database, the method greatly outperformed ANTS-SSD.
In LPBA40, IBSR18, CUMC12, and MGH10 databases, the method outperformed Diffeomorphic Demons.
In addition, the evaluation results in the regions of interest of Rohlfing et al. corroborated its excellent performance.

For Mang et al. benchmark PDE-constrained LDDMM methods~\cite{Mang_15,Mang_17b}, our evaluation results reported a significative loss of accuracy between 
RK and SL-RK integration. Interestingly, this loss of accuracy was not observed for our best performing method.

In IBSR18, CUMC12, and MGH10 databases, our PDE-constrained LDDMM methods were not able to reach SyN or Dartel performance.
This is because SSD image similarity metric restricts the performance of the methods in template-based segmentation.
This problem will be tackled in future work by formulating the PDE-constrained problem with other image similarity metrics such as
Normalized Cross Correlation, local Normalized Cross Correlation, or Mutual Information~\cite{Modersitzki_09_book}. 
We expect that this change in the formulation of the problem will increase the performance results of PDE-constrained LDDMM.
In addition, it will allow us to apply these methods to other clinical applications involving multi-modal registration. 

Simultaneously to the development of this work, Mang et al. released Claire software~\cite{Mang_19}. 
The software is intended to exploit massive CPU based parallel computing architectures to accelerate the computation time of PDE-LDDMM.
The codes implement original PDE-constrained LDDMM with a variational extension to nearly incompressible fluids and include $H^1$ and $H^2$ regularization terms.
The software is restricted to the stationary parameterization of diffeomorphisms.
The PDE integration scheme is SL-RK.
The software includes a sophisticated multi-level preconditioner that shows to improve the convergence of PCG with respect to the original proposal in~\cite{Mang_15}. 
The massive computation allows increasing the number of inner and outer iterations and use the norm of the gradient as stopping 
condition for achieving an extraordinary accuracy at convergence in a simulated experiment. 

In contrast to Claire, our methods of choice are intended to run completely in the VRAM of commodity GPUs ($<$ 4GBS).
We have limited the variational formulation to the one proposed in~\cite{Mang_15}, although it is straightforwardly extendible to the
nearly incompressible fluid problem.
We have limited our study to the traditional LDDMM regularizer.
Our software works for the stationary and the non-stationary parameterization of diffeomorphisms.
We have limited the preconditioner to the one proposed in~\cite{Mang_15} since we are interested on the comparison of the three different
variational variants.
We used the stopping conditions suggested in~\cite{Modersitzki_09_book} and used for PDE-LDDMM in~\cite{Mang_15,Mang_17,Hernandez_19,Hernandez_19b,Hernandez_19c}.
The variety of methods, the extensive evaluation conducted in this work, and our modest hardware capacity hindered us 
the use of the inner and outer iteration values needed for achieving the stopping conditions based on the norm of the gradient
suggested in~\cite{Mang_19}.
In fact, we observed in a selected NIREP experiment that increasing the number of inner iterations in PCG resulted into a faster initial convergence that finally 
stagnated in greater $MSE_{rel}$ values and lower $DSC$ scores than our considered stopping conditions.
This stagnation was also reported in~\cite{Mang_19} for the simulated experiment.
Instead, our selected inner and outer values consumed a reasonable amount of time while obtaining state of the art results for the 
evaluation metrics.

Comparing Claire and our results, we believe that it would be of interest to implement our best performing variant as a part of Claire's software.
In the other direction, it would be very interesting to adopt the multi-level preconditioners in our methods.

% De manera simultanea al desarrollo de este trabajo, Mang et al. liberaron los codigos de Claire.
% 
% - Variant 1 extendida a nearly incompressible
% - H1 results
% - Stationary parameterization only
% - Multi-level preconditioner que mejora el que se ha utilizado en Mang_15 y 16 y en Hernandez
% - Parametros de convergencia
% 
% ... De los resultados de este trabajo ...
% 
% - Merece la pena extender Claire a variants 2 and 3
% - Semi-Lagrangian integration in non stationary
% 
% - GPU vs CPU
% 
% - Parametros de convergencia

\section{Conclusions}
\label{sec:Conclusions}

In this work, we have proposed to combine the two different methodological approaches used to circumventing the huge computational complexity of 
Gauss--Newton--Krylov PDE-constrained LDDMM. 
In particular, we have included Semi-Lagrangian Runge-Kutta integration~\cite{Mang_17b} in the variants of band-limited PDE-constrained LDDMM proposed 
in~\cite{Hernandez_19,Hernandez_19c} 
for further increasing the computational efficiency of these methods. 
The resulting methods have been extensively evaluated in five different datasets following three different evaluation frameworks.
To our knowledge, this is the first time that SL-RK integration is implemented in the framework of PDE-LDDMM for the non-stationary parameterization and 
in the space of band-limited vector fields.
Moreover, our work first provides the position achieved by PDE-constrained LDDMM methods in the ranking of Klein et al. evaluation.

This study positions the formulation of BL PDE-constrained LDDMM based on the deformation state equation and SL-RK integration as the best performing
among all PDE-constrained LDDMM methods in terms of accuracy and efficiency.
The proposed method has reached the highest sensitivity in the classification of stable vs progressive mild cognitive impaired conversors in the Alzheimer's 
Disease Neuroimaging Initiative (ADNI) database using convolutional neural networks. This result has been recently published in~\cite{Ramon_20}. 

In future work, we will extend this formulation to other relevant physically meaningful LDDMM approaches such as the nearly incompressible method in~\cite{Mang_16}, 
and the geodesic shooting approach in~\cite{Hernandez_19b}. 
We will explore the advantages of using the multi-level preconditioner in~\cite{Mang_19}.  
We will adapt our methods for the use of alternative image similarity metrics that usually outperform SSD in registration evaluation rankings.
Finally, we will work in the understanding of which of the features of PDE-LDDMM allows the exceptional classification rates related with Alzheimer's disease
conversion shown in~\cite{Ramon_20}.

% \subsection{Evaluation cardiac (incompressible model!!)}

% - Comparacion de RK con SL: Monkey o paciente NIREP
% 
% - Evaluacion brain: NIREP, LPBA40, IBSR minimo
% 
% - Evaluacion cardiac: incompressible model (comparar con el diffeomorphic Demons?)

%-------------------------------------------------------

\section*{Acknowledgements}

The author would like to acknowledge the anonymous reviewers for their revision of the manuscript. 
% for their valuable comments and suggestions to improve the quality of the paper.
The author would like to give special thanks to Wen Mei Hwu from the University of Illinois 
for interesting ideas in the GPU implementation of the methods, 
and Nacho Navarro and Rosa Badia from the Barcelona Supercomputing Center (BSC) for his help.
%, for providing access and support to their multi-GPU cluster.
This work was partially supported by the national research grant TIN2016-80347-R.
In addition, this work was supported by NVIDIA through the Polytechnical University of Catalonia / 
Barcelona Supercomputing Center (UPC/BSC) GPU Center of Excellence.

%-------------------------------------------------------

%\begin{acknowledgements}
%\end{acknowledgements}


\bibliographystyle{splncs}      % mathematics and physical sciences

\bibliography{abbsmall.bib,Diffeo.bib,OpticalFlow.bib,Miccai2015.bib}

% % Non-BibTeX users please use
% \begin{thebibliography}{}
% %
% % and use \bibitem to create references. Consult the Instructions
% % for authors for reference list style.
% % Format for Journal Reference
% Author, Article title, Journal, Volume, page numbers (year)
% Author, Book title, page numbers. Publisher, place (year)
% \end{thebibliography}

\end{document}